\documentclass[11pt]{article} 

\usepackage{latexsym}
\usepackage{amsmath}
\usepackage{amssymb}

\newcommand{\vs}[1]{\vspace{#1mm}}

\newcommand{\hs}[1]{\hspace*{#1mm}}

\begin{document}

\hs{10}\begin{minipage}[t]{13.2cm} 
{\bf\Large Boolos--style proofs of limitative theorems}

\vs{2.5}

$\mbox{{\bf Gy\"orgy Ser\'eny}}^*$

\vs{1}

{\small 
Department of Algebra, 
Budapest University of Technology and Economics, 

\vs{-1}

1111 Stoczek\,u.\,2. \ H\,\,\'{e}p.\,\,5.\,em., 
Budapest, Hungary }

\vs{4}

{\small {\bf Key words:}\,\,Limitative results of logic, 
Boolos's  incompleteness proof, Berry's paradox

\vs{0}

{\bf MSC (2000)\,:}\,\,\,03F30}

{\vs 4}

\end{minipage}

\hs{10}\begin{minipage}[t]{13.2cm}
{\small Boolos's proof of incompleteness 
is extended straightforwardly to yield simple
``diagonalization--free'' proofs 
of some classical limitative theorems of logic.}

\end{minipage}

{\vs 4}

\hs{3}In his famous paper announcing the incompleteness theorem, 
G\"odel remarked that, though his argument 
is analogous to the Richard and the Liar paradoxes, 
``Any epistemological antinomy could be used for a similar proof
of the existence of undecidable propositions.'' 
([7]\,Note\,14). It is interesting that, despite the fact that 
the soundness of arguments like G\"odel's one built on \linebreak 
self--reference (or diagonalization) was often questioned 
(of course, from a philosophical not a mathematical point of view), 
the first attempt to support G\"odel's claim and prove the 
theorem using another paradox (and hence without recourse 
to diagonalization) came only recently. In 1989, 
formalizing the Berry paradox 
consisting in the fact that {\it the least integer not nameable in fewer 
than nineteen syllables} has just now been named in eighteen 
syllables, G.\,Boolos proved the semantic version of the incompleteness
theorem to the effect that there are arithmetical sentences 
that are true but unprovable in Peano arithmetic (see [3]). The proof, 
as Boolos notes at the end of his paper, 
``unlike the usual one, does not involve 
{\it diagonalization}''. Not much later, in a letter, he adds   
``What strikes the author as of interest in the proof via 
Berry's paradox is\,[\ldots]\,that it provides a {\it different 
sort of reason} for the incompleteness\,\,[\ldots]\,'' (cf.\,[4]). 

\hs{3}Perhaps Boolos's proof was one of the factors that 
have inspired a wave of ``proving old results in a new 
way'' (see e.g.\,\,[1] and the references given there). 
Nevertheless, unlike the proof theoretical methods used in  
both G\"odel's original proof and Boolos's one, 
most of these new proofs apply sophisticated model theoretical 
methods that can hardly be considered ``finitistic''.  
On the other hand, Boolos's proof can straightforwardly 
be extended to yield simple proofs of 
some fundamental theorems that are 
related very closely to the incompleteness theorem and 
to each other. The two versions of G\"odel's first  
incompleteness theorem (the semantical and syntactical one  
describing respectively the relation between truth and provability 
and that of provability and refutability)  
together with their strengthening 
(the G\"odel--Rosser theorem), Church's theorem on 
the undecidability of provability, 
and Tarski's theorem on the undefinability 
of truth, in a sense,  constitute a complete circle of mutually 
related statements answering some basic questions on 
provability and truth. The close connection between these 
fundamental results is also witnessed by the fact that 
their standard proofs have essentially the same structure\,: 
they all can be derived from a general formal version of the Liar 
paradox, that is, they can be considered as different 
formal resolutions of this paradox (cf.\,[11]). Now, as we shall
show below, almost the same can be said if we replace the Liar 
paradox by Berry's one. Actually, without any 
essential modification, the idea underlying Boolos's proof 
of incompleteness can be used to provide 
``diagonalization--free'' proofs 
of all the basic limitative theorems mentioned above.

\hs{3}After fixing notation and giving the definition of 
basic notions, we first mimic Boolos's proof in a 
slightly more detailed form than that in which 
the original proof was given so that 
we can continue the proof in different directions, 
which is just what we shall do.

\hs{3}Let us first fix any one of the standard first order languages 
of arithmetic. By a formula (resp.\,\,\,sentence, term etc.) 
we mean a formula (resp.\,\,\,sentence, term etc.) of this language. 
Theories are arbitrary sets of sentences. Robinson arithmetic 
(cf.\,[8]\,I.1.1)  will be denoted by $Q$. We shall 
denote the standard model of $Q$ (as well as its  
universe) by $\omega$, and say that a sentence is 
true (resp.\,\,\,a set is definable, defined etc.) 
if the sentence considered is true \vs{-2.5}
(resp.\,\,\,the set
is \linebreak 

\vs{-3}

\rule[0mm]{7cm}{0.1mm}
\vs{-1}

{\small ${}^*$ email\,: {\tt sereny\hs{0.2}@\hs{0.3}math.bme.hu}}

\pagebreak definable, defined etc.) in $\omega$.
The variables are $v_0,v_1,\ldots,v_i,\ldots\,\,$. 
If it seems necessary to indicate the difference between the 
closed terms \,$0, \,s\hs{0.5}0, \,s\hs{0.5}s\hs{0.5}0,\,\ldots \,$
and their values in $\omega$ (i.e. the natural numbers
0,1,2, \ldots), 
we shall denote the terms by the underlined versions of their
values, but as a rule, since there is no danger of confusion, 
we omit the underlining. Generally, the value of a closed 
term in $\omega$ will be 
denoted by the boldface version of the letter denoting the term 
concerned. 
Let us choose one of the standard G\"odel numberings. 
For any formula $\mu$, $\ulcorner\!\mu\!\urcorner$ 
will denote the G\"odel number of $\mu$. 
``iff'' stands for ``if and only if'', and we often use the symbol 
\,``\,\,$\circeq$\,\,''\, to stress that the equality 
concerned is a definition.

\hs{3}We say that a formula is $\Sigma_1$ if it is of the form 
$(\exists v_i)\,\mu$ for some $\Delta_0$ (i.e. bounded) formula $\mu$.
$\Sigma_1$ relations are those definable by a $\Sigma_1$ formula.   
A formula is $\Sigma$ (or a $\Sigma$ formula) 
if it is provably equivalent in $Q$ to some element of 
the smallest set (i.e. the intersection of all sets) containing  
all $\Delta_0$ formulas and being closed under conjunction, 
disjunction, existential quantification, and  bounded 
universal quantification. Further, a formula 
is said to be $\Delta$ (or a $\Delta$ formula)   
if both the formula itself and its negation are 
$\Sigma$\hs{0.2}. A $\Sigma$ ($\Delta$) sentence is a 
$\Sigma$ ($\Delta$) formula that is a sentence. 
Clearly, a $\Sigma_1$ formula is also a $\Sigma$ formula. 
A relation is called to be $\Sigma$ ($\Delta$) 
(or a $\Sigma$ ($\Delta$) relation) if it is definable 
by a  $\Sigma$ ($\Delta$) formula. It can easily be 
checked that a relation 
is $\Sigma$ iff it is $\Sigma_1$ 
(i.e. recursively enumerable), 
and a relation is $\Delta$ iff it is $\Delta_1$ 
(i.e. recursive), cf.\,e.g.\,[12]\,p.\,10. A straightforward 
induction on the complexity of formulas analogous to those 
that can be found in [5]\,(p.\,25) and [8]\,(I.1.8) shows  
that $Q$ is  $\Sigma$ complete, that is, 
all true $\Sigma$ sentences are provable in $Q$.\linebreak 

\vs{-2}

{\bf Definition}

\vs{2}

(i) For any term or formula $e$, let us denote by 
$|e|$ the number of symbols occurring in $e$ 
(we shall call this number the {\it length} of $e$),  
and let $f:\omega^2\longrightarrow \omega$ be a recursive function 
such that, for any formula $\mu$ and natural number $i$, 
$f(i, \ulcorner \!\mu\!\urcorner)
= \mbox{\raisebox{0.5mm}{$\ulcorner$}}\!
(\forall v_0)(\mu\!\! \iff \!\!v_0=i)\! 
 \mbox{\raisebox{0.5mm}{$\urcorner$}}$. 
(Obviously, there exists such a function.)

\vs{3}

(ii) For any theory $S$, let us denote by ${\mathcal Pr}_{\!S}$  
the set of G\"odel numbers of sentences provable in $S$. 

\vs{3}

(iii) Let $T$ be an arbitrary theory. 
Let us define the relations ${\mathcal Fm}\subseteq \omega$, 
and  ${\mathcal Lh},\,{\mathcal Nm}\subseteq \omega^2$
as follows:

\vs{3}

\hs{5}\begin{minipage}[t]{15cm}
${\mathcal Fm}\,\circeq\, \{i\in \omega: 
\,i= \ulcorner \!\mu\!\urcorner\mbox{ for some 
formula $\mu$ with at most one free variable 
$v_0$}\}$\,, 

\vs{2}

${\mathcal Lh}\,\,\circeq\,\, \{(i,j)\in \omega^2: 
\,i= \ulcorner \!\mu\!\urcorner\mbox{ for some formula $\mu$ 
such that  $|\mu|< j$}\}$\,,  

\vs{2}

${\mathcal Nm}\circeq \{(i,j)\in \omega^2: 
\,j\in {\mathcal Fm}\mbox{ and } f(i,j) \in {\mathcal Pr}_T\}$.  

\end{minipage}

\vs{3}

For any formula  $\mu$ and number $i$, 
if $(i, \ulcorner\!\mu\!\urcorner) \in 
{\mathcal Nm}$, that is, if $\mu\!=\!\mu(v_0)$ 
has at most one free variable 
$v_0$ and $T\vdash (\forall v_0)(\mu(v_0) \!\!\iff \!\!v_0=i)$, 
then we say that the formula $\mu$ {\it names} the number $i$.\linebreak

\vs{-2}

(iv) It follows from the definition of G\"odel numbering that 
the G\"odel numbers of formulas whose variables are all among the 
first ones are bounded by a recursive function of their length. 
More precisely, in the case of any G\"odel numbering, 
there is a recursive function $g$ (depending on the 
particular G\"odel numbering that has been chosen) such that,  
for any formula $\mu$ and number $j$,\linebreak whenever 
all the variables of  $\mu$ are among the first $j$ 
ones (that is, for $j\geq 1$, they are all in the set 
$\{v_0, v_1, \ldots, v_{j-1}\}$), $|\mu|< j$ implies that 
$\ulcorner\hs{-0.8}\mu\!\urcorner< g(j)$.\,\,\footnote{For example, 
let us consider the most commonly used G\"odel numbering,  
which (assuming that the G\"odel numbers of primitive 
symbols of the language concerned 
have already been given) is defined for any 
sequence of symbols as follows: 
$\mbox{\raisebox{0.6mm}{$\ulcorner$}}
\!\langle s_0, s_1, \ldots, s_j\rangle\!
\mbox{\raisebox{0.6mm}{$\urcorner$}}= 
p_0^{\ulcorner\! s_0\!\urcorner}\cdot 
p_1^{\ulcorner \!s_1\!\urcorner}\cdots 
p_j^{\ulcorner\!s_j\!\urcorner} $, where $p_i$ is
 the $i$th prime (see e.g.\,[10]\,pp.135--6). 
Apart from variables, our language has 
only finitely many primitive symbols, so we can define 
$c$ to be any number that is greater than the 
G\"odel numbers of primitive symbols except 
variables. \!Let $h(j)\circeq {\rm max}\,\{c\}\cup 
\{\ulcorner \!v_i\!\urcorner: i\!\leq\! j\}$ for every $j$.\linebreak   
Then, for any formula $\mu$ such that $|\mu|\!<\! j$ and 
all the variables of  $\mu$ are among the first $j$ ones, 
\,$\ulcorner \!\mu\!\urcorner < p_j^{\hs{0.3}h(j)\cdot j}$. 
Clearly,  the function $g(j)=p_j^{\hs{0.3}h(j)\cdot j}$ 
is \vs{-5}recursive.} Now, let us choose such a $g$, and 
let the relation ${\mathcal B}\subseteq \omega^2$ be defined 
in the following way: 

\vs{2}

\hs{3}${\mathcal B}\circeq 
\{(i,j)\in \omega^2: 
(\ulcorner \!\mu\! \urcorner, j) \in {\mathcal Lh} \mbox{ and }
(i, \ulcorner\hs{-0.2}\!\mu\! \urcorner) \in {\mathcal Nm}\, 
\mbox{ for some formula } \mu \mbox{ 
such that } \ulcorner \!\mu \!\urcorner <g(j)\}$.

\vspace*{2mm}

\pagebreak

(v) 
Obviously, ${\mathcal Fm}$ and ${\mathcal Lh}$ are $\Delta_1$ relations.
{\it Let us suppose that ${\mathcal Pr}_T$ is definable.} (This 
condition is obviously satisfied if, e.g., the set of G\"odel numbers of 
sentences belonging to $T$ is itself definable.) 
Then ${\mathcal Nm}$ and ${\mathcal B}$ are also definable.  

\vs{1.5}

\hs{5}(a) Let $\varphi(v_0,v_1)$ be a formula (with at most the  
free variables $v_0, v_1$) defining the relation 
${\mathcal B}$.\linebreak 

\vs{-3.5}

We shall choose $\varphi$ to be $\Sigma_1$ whenever 
$T$ is recursively axiomatizable. This is possible since, 
in this case, ${\mathcal Pr}_T$ is $\Sigma_1$, 
thus both ${\mathcal Nm}$ and 
${\mathcal B}$ are also $\Sigma_1$. 
(Recall that the class of recursively enumerable relations is closed 
under intersection, existential quantification, and the substitution 
of recursive functions, cf.\,e.g.\,[12]\,pp.\,27--8). 

\vs{2}

Note that, if a formula $\mu$ has at most one free variable $v_0$ 
and $|\mu|\!<\! j$, then, by renaming the bound variables of $\mu$,
we can obtain a formula $\mu^*$ such that $\mu^*$ 
has at most one free variable $v_0$, 
$\mu$ and $\mu^*$ \linebreak are provably 
equivalent in $Q$, $|\mu^*|=|\mu|$, and all the variables of 
$\mu^*$ are among the first $j$ ones, that \linebreak 
is, according to our 
remarks above, \raisebox{0.5mm}{$\ulcorner$}$\!\mu^*\!
$\raisebox{0.5mm}{$\urcorner$}$< g(j)$. 
In view of this fact, for any number $i$ and closed term $s$, 
\,$\varphi(i,s)$ is true iff  there is a formula $\mu$ 
such that $|\mu|< \mbox{{\bf s}}$ and $\mu$ names the 
number $i$. 

\vs{1.5}

\hs{5}(b) Let $\psi(v_0,v_1) \circeq \lnot\varphi(v_0,v_1)
\land(\forall v_2< v_0)\varphi(v_2,v_1)$. 

\vs{1.5}

For any number $i$ and closed term $s$, 
$\psi(i,s)$ is true iff  $i$ is the least natural number that 
cannot be named by a formula of length  $< \mbox{{\bf s}}$. 
(Clearly, $\psi(v_0,s)$ has at most one free variable $v_0$.)
 
\vs{3}

(vi)\,\, Let $k_1\circeq |\psi(v_0,v_1)|$ and let 
$k_2$ be any natural number that is greater than the number of 
free occurrences of $v_1$ in $\psi(v_0,v_1)$. 
Let $k\circeq k_1\cdot k_2$, \ 
$ t\circeq \underline{10}\cdot (\underline{k}\cdot \underline{k})$. \ 
Then \ $ k\geq k_1> 3, k\geq k_2 \geq 1$. 

\vs{2}

(vii) {\it If $T$ is a consistent extension of \,$Q$,} 
then every formula can name 
at most one number. (Indeed, $i\neq j$ implies 
$Q\vdash i\neq j$, cf.\,[8]\,I.1.6(3).)
Further, clearly, formulas provably equivalent in $T$ 
name the same number (if they name a number at all). 
Finally, up to provable equivalence 
in $T$, there are only finitely 
many formulas of less than a given length having at most one free 
variable $v_0$. (Recall that, apart from variables, our language has 
only finitely many primitive symbols and 
see\,\,our remarks in\,\,(v)\,(a).) 
Consequently, there are only finitely many 
different numbers that can be named by 
formulas of less than a given length. Thus, 
there is a least number that cannot be named by a formula of 
length less than ${\bf t}$. Let it be denoted by $n$.

\vs{3}

{\bf Theorem} 

\vs{1.5}

{\it If $T$ is a consistent extension of \,$Q$ and ${\mathcal Pr}_T$ 
is definable, then $\psi(n,t)$ is true, but $T\not\vdash\psi(n,t)$. } 

\vs{2.5}

{\large P\hs{0.3}r\hs{0.3}o\hs{0.3}o\hs{0.3}f.}
By definition, $n$ is the least number that cannot 
be named by a formula 
of length $< {\bf t}$,\linebreak  and, again by definition, 
$\psi(n,t)$ is true just in this case. Consequently,  

\vs{1}

(1) \ $\psi(n,t)$ is true.

\vs{1}

On the other hand, by the definition of $\psi$, (1) implies that 

\vs{1}

(2) \ $\varphi(n,t)$ is false.

\vs{1}

Further, it is easy to see that 

\vs{1}

(3) \ if \ $T\vdash \psi(n,t)$, \,then  
$\psi(v_0,t)$ names the number $n$. 

\vs{1}

Actually, we have to show that $T\vdash \psi(n,t)$ implies 
$T\vdash (\forall v_0)(\psi(v_0,t)\!\!\iff\! \!v_0=n$). 
Clearly, in one direction, the formal implication is trivial: 
$T\vdash \psi(n,t)\land v_0\!=\!n 
\Longrightarrow \psi(v_0,t)$. The other direction, in turn, 
follows from the fact that $T$ is an extension of $Q$ since 
$Q\vdash v_0\leq i\lor i\leq v_0$ (cf.\,[8]\,I.1.6\,(5)), which, 
in turn, implies a weak kind of provable uniqueness of least elements; 
more precisely, for {\it any} formula $\mu(v_0)$ and number $i$, 

\vs{0.5}

\ \hfill $Q\vdash  \lnot\mu(i)\land(\forall v_2< i)\mu(v_2)
\Longrightarrow (\forall v_0) [\lnot\mu(v_0)\land
(\forall v_2< v_0)\mu(v_2)
\Longrightarrow v_0=i\,]$. \hfill \ 

\vs{1}

Now, it follows from $k > 3$ that 
$18\hs{0.3}  k< 8 \hs{0.3} k^2$. Thus the definition of 
$t$ implies that $|\psi(v_0,t)|\leq |\psi(v_0,v_1)| + k_2\hs{0.3}|t|=
k_1+k_2(15+k+1+k+1)=k_1+k_2(17+2k)
\leq k+ k(17+ 2k)=
18\hs{0.3} k + 2\hs{0.3}k^2< 10\hs{0.3}  k^2={\bf t}$. 
So we have 
$|\psi(v_0,t)|< {\bf t}$, which, together with (3),  
shows that, if 
$T\vdash \psi(n,t)$, then $\psi(v_0,t)$ is actually 
a formula witnessing the truth of $\varphi(n,t)$. But 
$\varphi(n,t)$ is false by\,\,(2). Consequently, 
$T\not\vdash\psi(n,t)$.

\pagebreak

\hs{3}Now we can give the semantical incompleteness theorem in 
the usual formulation. The theory $T$ is called to be 
{\it sound} if all the sentences belonging to  $T$ are true.

\vs{4}

{\bf Corollary 1}  (Semantic version of G\"odel's first   
incompleteness theorem)

\vs{2}

{\it Let ${\mathcal Pr}_T$ be definable {\rm (}in particular,  
let $T$ be recursively axiomatizable{\rm )}. 
If $T$ is sound, then $T$ is incomplete.}

\vs{3}

{\large P\hs{0.3}r\hs{0.3}o\hs{0.3}o\hs{0.3}f.}
First of all, for any theory $S$, let us denote by 
Ded\,$S$ the set of all sentences 
provable \linebreak  in $S$, and let 
$T\,'\circeq Q\hs{0.2}\cup\hs{0.2}T$. 
Then obviously, Ded\,$T\,'$= Ded\,($Q\hs{0.2}\cup$\hs{0.2}Ded\,$T$), 
that is,  ${\mathcal Pr}_{T\,'}= 
{\mathcal Pr}_{Q\hs{0.2}\cup\hs{0.2}\mbox{{\scriptsize Ded}\,}T}$.
\linebreak  Since $Q$ is finite and ${\mathcal Pr}_T$ is 
definable by our assumption, the set of G\"odel numbers  
of the sentences in $Q\,\cup$\,Ded\,$T$ is again definable,
which, in turn, implies the definability of the set \linebreak  
\,${\mathcal Pr}_{Q\,\cup\,\mbox{{\scriptsize Ded}\,}T}= 
{\mathcal Pr}_{T\,'}$. Further,  
$Q$ is sound by definition, thus $T\,'$ 
is also sound. Soundness, in turn, implies consistency. 
Consequently, $T\,'$ satisfies the conditions of the 
Theorem. Therefore $T\,'\not\vdash\psi(n,t)$ 
and $\lnot\psi(n,t)$  is false. Thus, on the one hand, 
$T\not\vdash\psi(n,t)$ follows from the fact that $T\subseteq T\,'$, 
on the other hand, sound theories cannot prove 
false sentences. 

\vs{4}

\hs{3}So far we have only reiterated Boolos's proof 
with some minor modifications that open up the possibility 
to make a few steps farther along the lines set by the 
original proof, and  have formulated its most immediate 
consequence.\footnote{The detailed exposition, however, 
has its reward. The Theorem is a slightly 
more general version of the {\it semantical} incompleteness theorem 
than the usual one. Indeed, it seems that the standard proofs
(cf.\,e.g.\,[6]\,p.229 or [2]\,p.100), being 
essentially based on the diagonal lemma in one way or other, 
yield the theorem in such a form in which the condition of 
soundness of the theory concerned (which is, of course, a much 
stronger requirement than that of its consistency) inevitably appears; 
see the proof of the abstract version of this theorem in [11].  
For that matter, if we had followed Boolos's proof word by word, 
then we could have weakened even the condition that 
$T$ is an extension of $Q$. Actually, in order to define $n$, 
it is enough to suppose that the sentences 
$\underline{i}\neq \underline{j}, \,i,j \!\in\! \omega$ 
are all theorems of $T$. Then we can proceed as follows. 
$\psi(v_0,t)$ defines $n$ as the least element of a non--empty set of 
natural numbers. Since the least element of such a set is unique, 
the sentence $\eta\circeq (\forall v_0)
(\psi(v_0,t) \!\!\iff \!\!v_0=n)$, expressing the  uniqueness 
of this element and the fact that this element is just $n$, is true. 
On the other hand, the provability of 
$\eta$ in $T$ would entail that $\varphi(n,t)$ is true (recall that 
$|\psi(v_0,t)|<{\bf t}$), contradicting the truth of $\psi(n,t)$. 
As far as the third condition of the Theorem is 
concerned, the usual strong assumption of recursive axiomatizability 
of $T$ can obviously be weakened to the definability 
of ${\mathcal Pr}_T$ in the classical proofs as well.} 
In order to proceed, let us observe that, though
Boolos's proof is essentially a formalization of the Berry paradox, 
it is not the most straightforward one. As a matter of fact, the 
theorem that can be considered as the most faithful 
formal version of the Berry paradox is Tarski's theorem on the 
undefinability of truth. Rephrasing G\"odel's above quoted remark, 
we may conjecture that ``The formal version of any 
epistemological antinomy is just the statement on the 
undefinability of truth, and hence could be 
used for its proof''.  The reason is simple enough. 
As Tarski puts it in connection with the Liar paradox 
(cf.\,[13]\,p.\,76.), we cannot talk about the truth 
in the language of arithmetic 
since otherwise ``the antinomy of the liar could actually be 
reconstructed in this language''. As a simple corollary of the 
Theorem shows, literally the same can be said about 
the Berry paradox.

\vs{5}

{\bf Corollary 2} (Tarski's theorem on the undefinability 
of arithmetical truth)

\vs{2}

{\it The set of G\"odel numbers of true sentences 
is not definable}.

\vs{3}

{\large P\hs{0.3}r\hs{0.3}o\hs{0.3}o\hs{0.3}f.}
\,Let ${\mathcal T\!r}$ be the set of G\"odel numbers of 
true sentences and let us suppose that 
${\mathcal T\!r}$ {\it is} \linebreak definable. 
Choose $T$ in the Theorem to be the 
set of all true sentences, that is, let \linebreak
$T=\{\sigma: \sigma \mbox{ is a true sentence}\}$. 
Then, clearly, $T$ is deductively closed, i.e. the sentences provable 
in $T$ are all in $T$. Consequently, ${\mathcal Pr}_T={\mathcal T\!r}$, 
so that $T$ is consistent and ${\mathcal Pr}_T$ is definable.
Moreover, by definition, $T\supseteq Q$. So we can apply the 
Theorem: $\psi(n,t)$ is true but unprovable in $T$.
But ${\mathcal Pr}_T={\mathcal T\!r}$ implies that 
this is impossible because it means that, for any sentence 
$\sigma$, $\sigma$ is true iff $\sigma$ 
is provable in $T$.

\vs{4}

\pagebreak

\hs{3}Kikuchi has modified Boolos's 
notion of naming to obtain the syntactic version of the 
first incompleteness theorem for suitable extensions of Peano arithmetic 
and the second incompleteness theorem 
(see [9]).\footnote{There is, however, a 
minor mistake in his proof of the first incompleteness 
theorem (see the proof of Theorem 2.2\,(ii) in [9]). 
Indeed (using the notation of [9]), $Q(m,\rho)$ is obviously\,{\it not} 
$\Sigma_1$. What is needed, therefore, in order for 
{\it that} proof given in [9] to go through, is the simple fact 
(to be shown, needless to say, without using the soundness of $PA$) 
that there is a $\Sigma_1$ sentence $Q^*(m,\rho)$ satisfying not 
only the requirement that $Q(m,\rho)$ implies $Q^*(m,\rho)$ in $PA$, 
{\it but also} the additional one that the truth of $Q^*(m,\rho)$ 
implies the same for $Q(m,\rho)$.}
As far as the first incompleteness theorem is concerned, in fact, 
this modification is not needed. What is more important, 
with the help of making some plausible additional observations, 
we can derive the syntactic incompleteness theorem from the 
previous results for a considerably weaker theory than Peano arithmetic:

\vs{5}

{\bf Corollary 3} (Syntactic version of G\"odel's first 
incompleteness theorem)

\vs{2}

{\it Let $T$ be a recursively axiomatizable extension of \,$Q$. 

\vs{2}

{\rm (i)} \,\,If $T$ is consistent, then $T\not\vdash
\lnot\varphi(n,t)$. 

\vs{1}

{\rm (ii)} If $T$ is \,$\omega$--consistent, then 
$T\not\vdash\varphi(n,t)$. }

\vs{4}

{\large P\hs{0.3}r\hs{0.3}o\hs{0.3}o\hs{0.3}f.}
Since $T$ is recursively axiomatizable, $\varphi(v_0,v_1)$ is, 
by definition, a $\Sigma_1$ formula. Further, since 
$\omega$--consistency implies consistency, the conditions of 
the Theorem hold in both cases.

\vs{2}

(i) \,\,Since $\varphi(v_2,t)$ is  $\Sigma$\,, the sentence  
$(\forall v_2< n)\varphi(v_2,t)$ is also 
$\Sigma$. Moreover, it is true. (This follows from the definition of 
$\psi$ and (1) in the proof of the Theorem.) 
Therefore, by $\Sigma$ completeness, 
$T\vdash (\forall v_2< n)\varphi(v_2,t)$.
Hence $T\vdash \lnot\varphi(n,t)$ would imply 
$T\vdash \psi(n,t)$, contradicting the Theorem. 
Thus $T\not\vdash \lnot\varphi(n,t)$.\,\,\footnote{Note that 
$\lnot\,\varphi(n,t)$\, is yet another sentence that is true but 
unprovable.}

\vs{2}

(ii) \,Suppose that $T$ is $\omega$--consistent.  
Since $\varphi(v_0,v_1)$ is now $\Sigma_1$, 
there is a $\Delta_0$ formula $\mu(v_0, v_1, v_i)$ 
such that $\varphi(v_0,v_1)\!=\!(\exists v_i)\mu(v_0, v_1, v_i)$.
As we have already seen, $\varphi(n,t)$ is false
(cf. (2) in the proof of the Theorem). 
It follows from this that, for any number $j$,
$\lnot\mu(n,t,j)$ is a true $\Delta_0$ sentence. Using $\Sigma$ 
completeness, we have $T\vdash \lnot\mu(n,t,j)$ for every $j$, which,   
by the definition of $\omega$--consistency,\footnote{$T$ is 
{\it $\omega$--consistent} if for any formula 
$\eta(v_i)$, it follows from $T\vdash\!(\exists v_i)\eta(v_i)$ that 
$T\not\vdash\!\lnot\,\eta(j)$ for some number $j$ 
(cf.\,e.g.\,[10]\,p.\,142).}
implies that  $T\not\vdash (\exists v_i)\mu(n,t,v_i)$, i.e. 
$T\not\vdash\varphi(n,t)$. 

\vs{5}

\hs{3}One of the standard ways to prove the G\"odel--Rosser 
theorem is to show that it is a direct consequence of the 
Church theorem. We shall also follow this route, that 
is, using the previous results, we first show that no 
consistent extension of $Q$ is decidable\,: 

\vs{5}

{\bf Corollary 4} (Church's theorem on the undecidability 
of arithmetic)

\vs{2}

{\it If $T$ is a consistent extension of \,$Q$, then $T$ 
is undecidable.}

\vs{4}

{\large P\hs{0.3}r\hs{0.3}o\hs{0.3}o\hs{0.3}f.}
Suppose, for sake of contradiction, that 
${\mathcal Pr}_T$ is a recursive relation. 
It follows from this that ${\mathcal Nm}$ is also 
a recursive one, hence ${\mathcal B}$ is again recursive 
since, on the one hand, ${\mathcal Fm}$ and ${\mathcal Lh}$ 
are recursive, on the other, the class of recursive relations is closed 
under intersection, bounded quantification, and the substitution 
of recursive functions (cf.\,e.g.\,[12]\,pp.\,27--8). 
Consequently, the formula $\varphi$ defining ${\mathcal B}$ 
can now be chosen to be $\Delta$. Since $T$ 
is supposed to be a consistent extension of $\,Q$ 
and ${\mathcal Pr}_T $ to be recursive (which, of course, 
implies the recursive axiomatizability  
of $T$), \linebreak we can apply Corollary\,3\,(i). Consequently,  
$T\not\vdash\lnot\varphi(n,t)$. This, however, leads to a 
contradiction since, by (2) in the proof of the Theorem,  
$\lnot\varphi(n,t)$ is true, i.e. it is a true $\Sigma$ sentence. 
Its truth, in turn, by 
$\Sigma$ completeness, implies its provability in 
$T$, that is, $T\vdash \lnot\varphi(n,t)$. 

\vs{4}

\pagebreak 

\hs{3}In the usual way, Church's Theorem immediately yields 

\vs{4}

{\bf Corollary 5} (Rosser--G\"odel incompleteness theorem)

\vs{2}

{\it If $T$ is a consistent and recursively axiomatizable extension 
of \,$Q$, then $T$ is incomplete.}

\vs{4}

{\large P\hs{0.3}r\hs{0.3}o\hs{0.3}o\hs{0.3}f.}
Let us suppose that, on the contrary, $T$ is complete. 
Let ${\sf Snt}(v_0)$ and ${\sf Neg}(v_0,v_1)$ 
denote $\Delta$ formulas defining, respectively,  
the set of G\"odel numbers of sentences and the relation 
that holds between the G\"odel number of a sentence and that 
of its negation. Further, let ${\sf Pr}_T(v_0)$ a $\Sigma$ formula 
defining the set ${\mathcal Pr}_T$.  Now we set

\vs{1}
  
\ \hfill ${\sf Pr\hs{0.1}c\hs{0.3}}_T(v_0)
\circeq\lnot\,{\sf Snt}(v_0)\lor
(\exists v_1)({\sf Pr}_T(v_1)\land{\sf Neg}(v_0,v_1))$. 
\hfill \ 

\vs{1}

Then ${\sf Pr\hs{0.1}c\hs{0.3}}_T(v_0)$ is a $\Sigma$ formula. 
On the other hand, it follows from the completeness and consistency 
of $T$ that ${\sf Pr\hs{0.1}c\hs{0.3}}_T(v_0)$ 
defines just the complement of ${\mathcal Pr}_T$. Therefore,  
both ${\mathcal Pr}_T$ and its complement 
are $\Sigma$, so that ${\mathcal Pr}_T$ is recursive, contradicting 
the previous corollary. 

\vs{5}

\hs{3}The proofs we have given demonstrate that the Boolos--style 
formalization of Berry's paradox is, in fact, a proof {\it schema}. 
Indeed, in order to obtain the proofs of {\sl G\"odel's semantical 
incompleteness theorem, Tarski's theorem, G\"odel's syntactical 
incompleteness theorem,} and {\sl Church's theorem,} we have simply 
applied the {\it common} conceptual framework given implicitly by 
Boolos's incompleteness proof to four kinds of formal theories 
of arithmetic, namely, to theories for which the set 
${\mathcal Pr}_T$ is, respectively, {\sl definable, the set of 
{\rm (}G\"odel numbers of\hs{1}{\rm )}
all true sentences, recursively enumerable,} and {\sl recursive}. 

\vs{7}

{\bf Acknowledgment} 

\vs{2}

The research  was supported by Hungarian NSF grants No. T43242, 
T30314, and T035192.

\vspace*{5mm}

{\bf References}

\vs{3}

[1]\hs{3}Z. Adamowicz, and T. Bigorajska, Existentially closed 
structures and G\"odel's second 
 \linebreak
\hs{7}incompleteness theorem. 
J. Symbolic Logic {\bf 66}, 349--356(2001). 

\vs{1}

[2]\hs{3}D.W. Barnes and J.M. Mack, 
An Algebraic Introduction to Mathematical Logic (Springer--
\hs{7}Verlag, New York, 1975).

\vs{1}

[3]\hs{3}G. Boolos, A new proof of the G\"odel incompleteness theorem. 
Notices Amer. Math. Soc. 
\linebreak
\hs{7}{\bf 36}, 388--390(1989).

\vs{1}

[4]\hs{3}G. Boolos, A letter from George Boolos. Notices Amer. Math. 
Soc. {\bf 36}, 676(1989).

\vs{1}

[5]\hs{3}G. Boolos, The Logic of Provability (Cambridge University 
Press, Cambridge, 1995). 

\vs{1}

[6]\hs{3}H. B. Enderton, A Mathematical Introduction to Logic
(Academic Press, New York, 1972).   

\vs{1}

[7]\hs{3}K. G\"{o}del, On Formally Undecidable Propositions of 
{\it Principia Mathematica} and Related 
\linebreak
\hs{7}systems I. In: G\"{o}del's Theorem in Focus 
(S. G. Shanker, ed., Routledge, London 1988). 

\vs{1}

[8] P. H\'ajek, and P. Pudl\'ak, Metamathematics of First--Order
Arithmetic (Springer, \linebreak\hs{7}Berlin 1993).    

\vs{1}

[9]\hs{3}M. Kikuchi, A note on Boolos' proof of the incompleteness 
theorem. Math. Logic Quarterly \linebreak\hs{7}{\bf 40}, 
528--532(1994). 

\vs{1}

[10]\hs{1}E. Mendelson, Introduction to Mathematical Logic 
(D. Van Nostrand Company, Princeton \linebreak\hs{7}1964). 

\vs{1}

[11]\hs{1}G. Ser\'eny, \,\,\,G\"odel,\,Tarski,\,Church,\,and\,the\,Liar.
   The  Bulletin of Symbolic Logic {\bf 9}, 
3--25(2003).

\vs{1}

[12]\hs{1.5}R. M. Smullyan, Recursion Theory for Metamathematics
(Oxford Univ. Press, New York, \linebreak\hs{7}1993).

\vs{1}

[13]\hs{1.5}A. Tarski, Truth and proof. Scientific American  
{\bf 220},\,No\,6,\,63--77(1969).

\end{document}